\documentclass[11pt]{article}
\usepackage{amsmath} 
\usepackage{amssymb} 
\usepackage{latexsym}
\usepackage{theorem} 

\newtheorem{theorem}{Theorem} 

{\theorembodyfont{\rmfamily}

}
{\theoremstyle{break}	
{\theorembodyfont{\rmfamily} 
				\newtheorem{example}[theorem]{Example}
				
}}
{\end{enumerate}}

\newenvironment{proof}[1]{\smallskip \noindent {\bf #1}}{\qed\smallskip}

\def\qed{\ifhmode\unskip\nobreak\fi\ifmmode\ifinner\else\hskip5pt\fi\fi
 \hfill\hbox{\hskip5pt\vrule width4pt height6pt depth1.5pt\hskip1pt}}

\newcommand{\ov}[1]{\ensuremath{\overline{#1}}}

\newcommand{\RR}{\ensuremath{{\bf R}}}  
\newcommand{\R}[1]{\ensuremath{{\bf R}^{#1}}} 

\newcommand{\NN}{\ensuremath{{\bf N}}}	   
  
\newcommand{\om}[1]{\ensuremath{\omega(#1)}}

\newcommand{\pde}[2]
	{\ensuremath{\frac{\partial #1}{\partial #2}}} 
\newcommand{\ode}[2]{\ensuremath{\frac{d#1}{d#2}}}    	

\newcommand{\abs}[1]{\ensuremath{\vert#1\vert}}
\def\d#1dt{\frac{d#1}{dt}}    

\newcommand{\lam}{\ensuremath{\lambda}}

\newcommand{\del}{\ensuremath{\delta}}

\newcommand{\sig}{\ensuremath{\sigma}}
\newcommand{\Sig}{\ensuremath{\Sigma}}
\newcommand{\Gam}{\ensuremath{\Gamma}}
\newcommand{\gam}{\ensuremath{\gamma}}
\newcommand{\co}{\colon\thinspace} 

\newcommand{\msf}[1]{{\mathsf {#1}}}
\newcommand{\mbf}[1]{{\mathbf{#1}}}

\newcommand{\dimply}{\ensuremath{\:\Longleftrightarrow\:}}

\def\mylabel#1{\label{#1}}
\reversemarginpar

\def\I{\mbf I}

\def\K{\msf K}
\def\Ko{\msf K^o}

\def\p{\partial}

\newcommand{\comment}[1]{}

\begin{document}

\title{ On existence and uniqueness of the carrying simplex for 
  competitive dynamical systems}
\author{Morris W. Hirsch\\
 University of California, Berkeley\\
University of Wisconsin, Madison
}
\date{}

\maketitle
 \comment{ \begin{flushright} {\footnotesize Npapers/CarryingSimplex/jbd3. tex}
  \end{flushright}
  \normalsize}


\begin{quote}{\em Dedicated to Professor Hal Smith on the occasion of
    his sixtieth birthday}
\end{quote}

g
\section*{Introduction}
Consider a system of $n$ competing species whose states
are characterized  by vectors in  the closed positive cone 
$\K:=[0,\infty)^n\subset \R n$.  When time is discrete the development
  of the system is given by a continuous map $T\co \K \to \K$.
For continuous  time is continuous the development is governed by a
periodic system of  differential equations $\dot x = F(t,x)\equiv
F(t+1,x)$.  In this case $T$ denotes the  Poincar\'e map.   

For discrete time the {\em trajectory} of a state $x$ is the sequence
$\{T^k x\}$, also denoted by $\{x(k)\}$, where  $k$ varies over the set
$\NN$ of nonnegative integers.  In the case of an autonomous
differential equation (i.e., $F$ is independent of $t$),  the
trajectory of $x$ is the solution curve through $x$, denoted by $T^tx$
or $x(t)$, where $t\in[0,\infty)$.  In both cases the
{\em limit set} $\om x$ is the set of limit points of sequences
$x(t_k)$ where $t_k\to \infty$.  

In order to exclude 
spontaneous generation we assume $T_i(x)=0$ when $x_i=0$.  
Thus there are 
functions $G_i\co \K\to [0,\infty)$, assumed continuous, such that
\begin{equation}		\label{eq:TG}
 T_i(x)=x_i G_i(x), \qquad ( x\in \K, \quad i=1,\dots,n)
\end{equation}
For continuous time we assume the differential equation is a system of
having the form $\dot x_i=x_iG_i(t,x)$.  
If $x_i$ is interpreted as the size of species $i$ then $G_i (x)$ is
its {\em per capita} growth rate.

We take ``competition'' to mean that  increasing any one species does
not tend to increase the {\em per capita} growth rate of any other
species, conventionally modeled by the assumption
\[
  \textstyle \frac{\p G_i}{\p x_j}\le 0, \qquad(i\ne j)
\]

A {\em carrying simplex} for the map $T$ is a compact invariant hypersurface
$\Sig\subset \K$ such that every trajectory except the the origin is
asymptotic with a trajectory in $\Sigma$, and $\Sigma$ is unordered
for the standard vector order in $\K$.  In the case of an autonomous
differential equation we require that $\Sig$ be invariant under the 
maps $T^t$ for all $t\ge 0$. 
Some maps have no carrying simplices, others have infinitely many.
Our main results gives conditions guaranteeing  a
unique carrying simplex.  

\subsection*{Terminology}
A set $Y\subset \K$ is {\em positively invariant} under a map or an
autonomous differential equation if it contains the
trajectories of all its points, so that 
$T^t Y\subset Y$
for all $t\ge 0$ (here $t\in \NN$ or $[0,\infty)$ as appropriate).  We
  call $Y$  {\em invariant} if  $T^t Y= Y$ for all $t\ge 0$. 

If  $S$ is  a differentiable map, its matrix of partial
derivatives matrix at $p$ is denoted
by $S'(p)$.  

The geometry of $\K$ plays an important role.  For each subset
$\mbf I\subset\{1,\dots,n\}$ the $\I$'th  {\em facet} of $\K$ is 
\[
  \K_\mbf I=\{x\in \K\co x_j =0 \dimply j\notin \mbf I\}
\]
Thus $\K_{\{i\}}$ is the $i$'th positive
coordinate axis. 
A facet is {\em proper} if it lies in the boundary of $\K$, meaning
$\I\ne \{0\}$.  
The intersection of facets is a facet:
 $\K_I\cap\K_J=\K_{I\cup J}$.    The boundary of
 $\K$ in $\R n$, denoted by $\dot \K$, 
is the union of the proper facets.
Each $x\in \K\verb=\=\{0\}$ belongs to the  unique facet 
    $\K_{\I(x)}$ where   $\I(x):=\{i\co x_i > 0\}$.  

For each  $n\times n$ matrix $A$ and nonempty
$\I\subset\{1,\dots,n\}$ we define the principal submatrix
\[
  A_\I:=\left[A_{ij}\right]_{i, j\in \I}
\]

The {\em vector order} in $\R n$ is the relation is defined by
$x\succeq y \dimply x-y \in \K$.  We write $x\succ y$ if also $x\ne
y$.  For each set $\I\subset\{1,\dots,n\}$ we write $x\succ_\mbf I
y$ if $x, y\in 
\ov{\K_\mbf I}$ and $x\succ y$, and $x\succ\succ_\mbf I y$ if also $ x_i>y_i$
for all $i\in\mbf
I$.  The reverse relations are denoted by $\preceq ,
\prec $ and so forth.  

The {\em closed order interval} defined by $a,
b\in \R n$ is
\[
 [a,b]:=\{x\in\R n \co a\preceq  x\preceq  b\}
\]

\section*{Carrying simplices} 
A {\em carrying simplex} is a set $\Sigma\subset
\K\verb=\=\{0\}$ having the following properties:
\begin{description}

\item  [(CS1)] $\Sigma$ is  compact and invariant.

\item  [(CS2)] for every $x\in \K\verb=\= \{0\}$ the trajectory of $x$
  is asymptotic with some $y\in
\Sigma$,\ i.e., $\lim_{t\to\infty}|T^t x-T^ty|=0$.

\item [(CS3)] $\Sig$ is {\em unordered}: if $x, y \in \Sig$ and
  $x\succeq y$ then $x=y$.  

It follows that each line in $\K$ through the origin
  meets $\Sig$ in a unique point.  Therefore  $\Sig$ is mapped
  homeomorphically onto the
  unit (n-1)-simplex 
\[ \Delta^{n-1}:=\{x\in \K\co\sum_i x_i=1\}\]
 by the  radial projection $x\mapsto x/(\sum_i x_i)$.

\end{description}

Long-term dynamical properties of trajectories are  accurately reflected 
by the dynamics in $\Sig$ by (CS1) and (CS2),  and  (CS3) means that $\Sigma$
has  simple topology and geometry.  
The existence of a carrying simplex has significant implications for 
limit sets $\om x$:
\begin{itemize}

\item If $x >0$ then $\om x \subset \Sig$, a consequence of (CS2).  In
  particular, $\Sig$ contains all nontrivial fixed points and periodic
  orbits.

\item If $a, b \in \K$ are distinct limit points of respective states
  $x, y \succeq 0$ (possibly the same state), then there exist $i, j$
  such that $a_i >b_i, \ a_j < b_j$; this follows from (CS3).  Thus
  either $\om x =\om y$, or else there exist $i, j$ such that 
\[
  \limsup_{t\to\infty}\, x_i(t)- y_i(t) > 0, \qquad
  \liminf_{t\to\infty}\,  x_j(t)- y_j(t)< 0
\]
\end{itemize}

In many cases $\Sig$ is the {\em
global attractor} for the dynamics in $\K\verb=\=\{0\}$, meaning that
as $t$ goes to 
infinity, the distance
from $x(t)$ to  $\Sig$ goes to zero uniformly for
$x$ in any given compact subset of $\K\verb=\=\{0\}$.  This implies
(Wilson \cite{Wilson69}) that there is a continuous function $V\co
\K\verb=\=\{0\}\to [0,\infty)$ such that if $x\ne 0$ then
\begin{itemize}
\item
$V(x)=0\dimply x\in\Sig$,

\item $V(x(t)) < V(x) \dimply x\notin \Sig$,

\item $\lim_{t\to\infty} V(x(t))=0$,
\end{itemize} 
We can think of $V$ as an ``asymptotic conservation
law''.  While there are many such functions for any carrying simplex,
it is rarely possible to find a formula for any of them.

Before stating results we give two simple examples for $n=1$: 
\begin{example}  \mylabel{th:ex00}
If $T$ is the time-one map for the flow defined by the logistic
differential equation
\[
 \dot x =rx (\sigma -x), \quad  r, \sigma>0,\qquad (x\ge 0),
\]
the carrying simplex is just the classical  carrying capacity
$\sigma$.   Here one can define $V (x)= |x-\sig|$ for $x>0$.    
\end{example}

\begin{example}  \mylabel{th:ex0}
Consider the map 
\begin{equation}		\label{eq:0}
T\co [0,\infty) \to [0, \infty),\ Tx= xe^{b-ax},  \quad  b, a >0,
  \qquad x\in [0,\infty) 
\end{equation}
Note that
\[
  T' (x) = (1-x) e^{b-ax},\qquad T' (b/a) = 1-b 
\]
If  there is a carrying simplex,    it has to be 
the unique positive fixed point $b/a$, in which case
$\lim_{k\to\infty} T^k x = b/a$ for all $x >0$.   
 
{\em If $b \le 1$ then $b/a$ is the carrying simplex. }   In this case the
maximum value of $T$ is taken uniquely at $1/a \ge b/a$.  If 
$0<x<b/a$ then $x < Tx <b/a$, hence $T^kx \to b/a$.  It follows that
if some $T^j x  <b/a$ then again $T^k x \to b/a$.  If the entire orbit
of $x$ is $> b/a$ then the sequence $\{T^k x\}$ decreases to a fixed
point $\ge b/a$, hence to be $b/a$.

{\em If $b>2$ there is no carrying simplex. }   For then
$|T'(b/a)| > 1$, making  $b/a$ a locally repelling fixed point. 
The only way
the trajectory of  $y \ne b/a$ 
can converge to $b/a$ is for $T^j y= b/a$ for some $j >0$.      
The set of such
points $y$ is nowhere dense because $T$ is a nonconstant analytic function,
hence there is no carrying simplex.   For sufficiently large $b$ the
dynamics is chaotic.

Example \ref{ex:exMay}, below,  is an $n$-dimensional generalization of
Equation (\ref{eq:0}). 
\end{example}

We say that $T$ is {\em strictly  sublinear}  in a set $X\subset \K$
if the following holds: $x\in X$ and $0<\lam<1$ imply $\lam x \in X$ and 
\begin{equation}		\label{eq:3bis}
 \lam T (x)  \prec T(\lam x), \quad (x\in X \setminus{0})
\end{equation}
Thus the restricted map $T|X$ exhibits what economists call
``decreasing returns to scale.''

A state $x$ {\em majorizes} a state $y$ if $x \succ y$, 
and $x$ {\em  strictly  majorizes} $y$ if $x_i >0$ implies $x_i >
y_i$. 

The map $T\co\K\to\K$ is  {\em strictly retrotone} in a subset
$X\subset \K$ if
for all $x, y \in X$ we have
\[
  \mbox{$Tx$ majorizes $Ty\ \implies \ x$  strictly majorizes $y$}
\]
Equivalently: 
\[
 \mbox{ $x, y\in X\cap \ov{\K_\mbf I}$ \ and \ $Tx\succ Ty \ \implies \
 x\gg_\mbf I y$}
\]

The origin is a {\em repellor} if $T^{-1}(0)=0$ and there exists $\del >0$
and an open neighborhood $W\subset\K$ of the origin such that
$\liminf_{k\to\infty} \abs {T^kx} \ge \del$ uniformly in compact
subsets of $W\verb=\=\{0\}$.   

If in addition there is a global
attractor $\Gam$, as will be generally assumed, then $\Gam$ contains a
global attractor $\Gamma_0$ for $T|\, \K\verb=\=\{0\}$.  In 

We will assume $T$ is given Equation (\ref{eq:TG}) has the following
properties:  
\begin{description}

\item[(C0)]  $T^{-1}(0)= 0$ and $G_i (0) > 1$.  

The first condition is means that no nontrivial population dies out in
finite time.  The second means that small populations increase.  

\item[(C1)]  {\em There is  a global attractor $\Gam$ containing a
  neighborhood of $0$. } 

Together with (C0) this implies that there is a global attractor
$\Gamma_0 \subset \Gam$ for $T|\, \K\verb=\=\{0\}$.  The connected
component of the origin in $\K\verb=\=\Gam_0$ is the {\em repulsion
basin} $B(0)$.

\item[(C2)] {\em $T$ is  strictly  sublinear in a neighborhod of $\Gam$.}   

This holds when $
  0<\lam <1\implies   G (x) \prec G (\lam x)$.

\item[(C3)] {\em $T$ is strictly retrotone in a
  neighborhood of the  global attractor} 

A similar property was introduced by Smith \cite{Smith86}.
\end{description}

Denote the set of
boundary points of $\Gamma$ in $\K$ by $\p_\K\Gam$.
\begin{theorem}		\mylabel{th:mainMAPS}
When {\em (C0)---(C3)} hold, the unique carrying simplex is
$\Sig =\p_\K\Gam=\p_\K B(0)$, and $\Sig$  is the global attractor 
for $T|\, \K\setminus\{0\}$.
\end{theorem}
The proof will appear elsewhere. 

The same hypotheses yield further information.  It turns out that if
$T |\Gam$ is locally injective (which Smith assumed), it is a
homeomorphism of $\Gam$;  and in any case the following condition
holds:

\begin{description}

\item[(C4)] {\em The restriction of $T$ to each positive coordinate
axis   $\Ko_{\{i\}}$   has a globally attracting fixed point $q_{(i)}$.}

\end{description}
We call $q_{(i)}$ an {\em axial} fixed point.  Denoting its $i$'th
coordinate by $q_i >0$, we set 
\[
  q:=(q_1,\dots,q_n)=\sum_i q_{(i)} 
\]
Smith \cite{Smith86} shows that (C3) and (C4) imply (C1) with
$\Gam\subset [0, q]$.  In many cases the easiest way to establish a
global attractor is to compute the axial fixed points and
apply Smith's  result. 

The following condition implies (C3) for  maps 
$T$ having the form (\ref{eq:TG}) when 
$G$ is $C^1$:
\begin{description}
\item[(C5)] {\em If $x\in \K\verb=\=\{0\}$, 
the matrix $\left[G'(x)\right]_{\mbf I (x)}$ has strictly
  negative entries} 
\end{description}

For $d\in \R n$ we denote the diagonal matrix $D$ with diagonal
entries $D_{ii}:=d_i$ by $[d]^{\msf {diag}}$ and also by $[d_i]^{\msf
{diag}}$.
The $n\times n$ identity matrix is denoted by $I$.

A  computation shows that 
\[
  T'(x)=[G (x)]^{\msf{diag}} + [x]^{\msf {diag}}G' (x).
\]
When $x$ is such that all $G_i (x) >0$, this can be written 
\begin{equation}		\label{eq:2b}
\begin{split} T'(x)&=[G(x)]^{\msf{diag}}(I- M(x)),\\
 M(x) &:=  -\left[\frac{x_i}{G_i (x)}\right]^{\msf{diag}}G' (x),
\end{split}
\end{equation}
and the entries in the $n\times n$ matrix $M(x)$ are
\begin{equation}		\label{eq:2a}
\begin{split}  
   M_{ij}(x)&:= 
\frac{-x_i}{G_i (x)}\frac{\p G_i}
    {\p x_j} (x),\\
&= 
  -x_i\frac{\p\log G_i}{\p x_j}(x)
\end{split}
\end{equation}
Note that (C5) implies $M_{ij}(x)> 0$. 

The {\em spectral radius} $\rho (M)$ of an $n\times n$ matrix $M$ is
the maximum of the norms of its eigenvalues.  It is a standard result
that if $\rho(M) <1$ then $I-M$ is invertible  and
$(I-M)^{-1}=\sum_{k=0}^\infty M^k$.  
\begin{theorem}		\mylabel{th:maps4}
Suppose $G$ is $C^1$. Assume {\em (C0), (C1), (C2), (C5)}, let {\em (C4)}
 hold with $[0, q]\subset X$, and assume
\begin{equation}		\label{eq:4}
0 \prec x \preceq  q \implies \rho (M(x)) <1
\end{equation} 
Then {\em (C3)} holds, whence  the hypotheses and conclusions of Theorem
\ref{th:mainMAPS} are valid.
\end{theorem}
The proof will be given elsewhere. 
Under the same hypotheses the following  conclusions also hold:
\begin{itemize}
\item
 {\em $T|\Gamma$ is a diffeomorphism}

\item  {\em if $x\in\Gamma \cap \Ko_\I$ then 
the matrix $[T'(x)_{\mbf I}]^{-1}$ has strictly positive entries.}

\end{itemize}

When (C5) holds,  either of the following conditions implies
(\ref{eq:4}):
\begin{equation}		\label{eq:3a}
0\prec x\preceq  q\implies \sum_i M_{ij}(x) <1, \quad (j=1,\dots,n)
\end{equation}
\begin{equation}		\label{eq:3b}
0\prec x\preceq  q\implies \sum_j M_{ij}(x) <1 \quad (i=1,\dots,n)
\end{equation}
Each of these conditions implies that the largest positive eigenvalue
of $M(x)$ is the spectral radius by (C5) and the theorem of Perron and
Frobenius \cite{BermanNeumann89}, and that this eigenvalue is bounded
above by the maximal row sum and the maximal column sum by
Gershgorin's theorem \cite{BrualdiMellendorf94}.

\section*{Competition models}
In the following illustrative examples we calculate bounds on
parameters that make row sums of $M (x)$ obey (\ref{eq:3a}),
validating the hypotheses and conclusion of Theorem \ref{th:maps4} and
\ref{th:mainMAPS}.
\begin{example}  \mylabel{ex:exMay}
Consider a multidimensional version of Equation (\ref{eq:0}),  based
on an ecological model of  May \&  Oster \cite{May}: 
\begin{equation}\label{eq:may}
T\co \K\to \K,\quad T_i (x) = x_i\ {\textstyle \exp\big(B_i-\sum_j
  A_{ij} x_j\big)}, \qquad
B_i, \;A_{ij}>0
\end{equation}
This map is not locally injective.  In a small neighborhood of the
origin $T$ is approximated by the discrete-time Lotka-Volterra map
$\hat T$ defined by $(\hat T x)_i = (\exp B_i) x_i (1-\sum_j
A_{ij}x_j)$, but as $\hat T$ does not map $\K$ into itself, it is not
useful as a global model.  $T$ has a global attractor $\Gamma$ and a
source at the origin, so a carrying simplex is plausible.  But the
special case $n=1$, treated in Example \ref{th:ex0}, shows that
further restrictions are needed.

Condition (C5)  holds  with
$
  G_i (x)= {\textstyle \exp\big(B_i-\sum_j
  A_{ij} x_j\big)}$.
Evidently these  functions  are strictly
decreasing in $x$, which implies $T$ is strictly sublinear.
 (C4) holds with $q_i=B_i/A_{ii}$, 
and it can be shown that $\Gamma \subset [0, q]$.
In (\ref{eq:2b}) the 
 matrix entries  are  
\begin{equation}		\label{eq:mij}
M_{ij}(x)=x_i A_{ij}
\end{equation}
Therefore  Theorem \ref{th:maps4} shows that if
\begin{equation}		\label{eq:rho}
0\prec x\preceq  q \implies \rho (M(x)) < 1
\end{equation}
then $\p_\K\Gam$ is the unique carrying simplex and $T|\Gamma$ is a
 diffeomorphism.  From (\ref{eq:3a}), (\ref{eq:3b}) and (\ref{eq:mij})
we see that (\ref{eq:rho}) holds in case one of the following
 conditions is satisfied:
\begin{equation}		\label{eq:rho2}
 \mbox{$\displaystyle \frac{B_i}{A_{ii}}\sum_j A_{ij} < 1$ for all
  $i$,}
\end{equation}
 or 
\begin{equation}		\label{eq:rho2a}
\mbox{$\displaystyle \frac{B_i}{A_{ii}}\sum_i A_{ij}  < 1$ for all
     $j$}
\end{equation}
These conditions thus imply a unique carrying simplex, by Theorem
\ref{th:maps4}.  

To arrive at a biological interpretation of 
(\ref{eq:rho2}), we rewrite it as 
\begin{equation}		\label{eq:rho3}
 q_i \sum_j A_{ij} <1 
\end{equation}
where $q_i :=\frac{B_i}{A_{ii}}$ is the axial equilibrium for species
$i$, that is, its stable population in the absence of competitors.
Equation (\ref{eq:may}) tells us  that  $A_{ij}$ is the logarithmic
rate by which the growth of population $i$ inhibits the growth rate of
population $j$.  Thus (\ref{eq:rho3}) means that  the average of these
 rates must be rather small compared to the single species
equilibrium for population $i$.  The  plausibility of
this x1is left to the reader, as is the biological
meaning of  (\ref{eq:rho2a}).
 
When $n=1$, Equation (\ref{eq:may}) defines the map $Tx= xe^{b-ax}$ of
Example \ref{th:ex0}.  The positive fixed point is $q=a/b$, and both
(\ref{eq:rho2}) and (\ref{eq:rho2a}) boil down to $b<1$, which was
shown to imply a unique carrying simplex.  That example also showed
that there is no carrying simplex when $b>2$.  As Equation
(\ref{eq:may}) reduces to Example \ref{th:ex0} on each coordinate
axis, we see that Equation (\ref{eq:may}) lacks a carrying simplex
provided
\[
   \mbox{$\displaystyle \frac{B_i}{A_{ii}}\sum_j A_{ij} >2$ for some 
  $i$,}
\]
or 
\[
  \mbox{$\displaystyle \frac{B_i}{A_{ii}}\sum_i A_{ij}  > 2$ for some
     $j$}
\]
\end{example}

\begin{example}  \mylabel{th:exgower}
Consider a competing population model due to Leslie \& Gower
\cite{LeslieGower}:

\begin{equation}		\label{eq:lesgow}
T\co\K\to\K,\quad T_ix = \frac{C_i x_i}{ 1+ \sum_j A_{ij}x_j},\qquad
C_i,\;A_{ij}>0 
\end{equation}
Note that $T$ need not be locally injective.   When $n=1$ all
trajectories converge to $0$ if $C \le 1$, and all nonconstant
trajectories converge to $\frac{C-1}{A}$ if $C >1$. 
The case $n=2$ is thoroughly
analyzed by Cushing {\em et al.\ }\cite{Cushing04}.

Here 
\[
  G_i (x):= \frac{C_i}{ 1+ \sum_j A_{ij}x_j} >0,
\]
hence (C5) holds. 
We assume $C_i>1$, guaranteeing (C4)  with $q_i=
  \frac{C_i-1}{A_{ii}}$.
In (\ref{eq:2a}) we have 
\[
  M_{ij}(x)= \frac{x_i A_{ij}}{1+\sum_l A_{il}x_l}
< x_i A_{ij}. 
\]
so the row sums of $M (x)$ are $<1$ for all $x$ provided $q_i\sum_j
A_{ij} <1$.  Therefore  when
\[
  1<C_i < 1+ \frac {A_{ii}}{\sum_j A_{ij}},
\] 
Theorems \ref{th:mainMAPS} and \ref{th:maps4} yield the following
conclusions:  There is a global
attractor $\Gam\subset [0,q]$, the unique carrying simplex is 
$\p_\K \Gam$, and $T|\Gam$ is a diffeomorphism.
\end{example} 

\begin{example}  \mylabel{th:exneural}
Consider a recurrent, fully  connected neural network of $n$ cells
(or ``cell assemblies'', Hebb \cite{Hebb49}).  At discrete times
$t=0,1,\dots$, cell $i$ has activation level $x_i(t)\ge 0$ and the
state of the system is $x(t):=(x_1(t),\dots,x_n(t))$.  Cell $i$
receives an input signal $s_i(x(t))$ which is a weighted sum of all
the activations plus a bias term.  Its activation is multiplied by a
positive transfer function $\tau_i$ evaluated on $s_i$, resulting in
the new activation $x_i (t+1)=x_i (t)\tau_i(s_i)$.

We assume each cell's activation tends to decrease the activations of
all cells, but each cell receives a bias that tends to increase its
activation.  We model this with negative weights $-A_{ij}<0$, positive
biases $B_i>0$, and positive increasing transfer functions.  For
simplicity we assume all the transfer functions are $e^\sig$ where
$\sig\co [0,\infty)\to [0,\infty)$ is $C^1$.  States evolve according
to the law
\[
  T\co \K\to\K, \quad T_i(x)=x_i\exp \sig (s_i(x)),\qquad
  s_i(x):= B_i-\sum_j   A_{ij}x_j
\]
We also assume
\begin{equation}		\label{eq:neural2}
\sig(0)=0,\quad \sig'(s) >0, \quad \sup \sig'(s)=\gam
<\infty, \qquad (s\in\RR)
\end{equation}
It is easy to verify that (C1), (C2), (C4) and (C5) hold,   with 
\begin{equation}		\label{eq:neural3}
q_i:= \frac{B_i}{A_{ii}},\quad  
G_i (x):=\exp \sig \big(B_i-\sum_j
A_{ij}x_j\big),\quad M_{ij}= \sig' (s_i(x))A_{ij}
\end{equation}
where $M_{ij}(x)$ is defined as in (\ref{eq:2a}). 

It turns out that for given weights and biases, the system has a unique
carrying simplex provided the {\em gain parameter} $\gamma$ 
in (\ref{eq:neural2}) is not too large.  It suffices to assume 
\begin{equation}		\label{eq:neural4}
 \gam < \left[\max_{i}\bigg(\frac{B_i}{A_{ii}}\sum_j
 A_{ij}\bigg)\right]^{-1} 
\end{equation}
For then (\ref{eq:neural2}), (\ref{eq:neural3}), (\ref{eq:neural4})
imply (\ref{eq:3b}) and hence (C3), so Theorems
\ref{th:mainMAPS} and \ref{th:maps4} imply a unique carrying simplex for
$T$. 

There is a vast literature on neural networks, going back to the
seminal book of Hebb \cite{Hebb49}.  Network models of competition
were analyzed in the pioneering works of Grossberg \cite{Grossberg78}
and Cohen \& Grossberg \cite{CohenGrossberg83}.  Generic
convergence in certain types of competitive and cooperative networks
is proved in  Hirsch
\cite{Hirsch89}.  Levine's book \cite {Levine00} has mathematical
treatments of several aspects of neural network dynamics.
\end{example}

\section*{Competitive differential equations}
Consider a periodic differential equation in $\K$:
\begin{equation}		\label{eq:per}
\dot u_i = u_iG_i(t, u_1,\dots,u_n)\equiv u_iG_i(t+1, u_1,\dots,u_n),
 \quad t, u_i \ge 0, \quad (i=1,\dots,n))
\end{equation}
where the maps $G_i\co \K\to\RR$ are $C^1$.
The solution with initial
value $u(0)=x$ is denoted by $t\mapsto T_t x$.  Solutions are assumed
to be defined for all $t \ge 
0$.  Each map $T_t$ maps $\K$ diffeomorphically onto a relatively open
set in $\K$ that contains  the origin.   The {\em Poincar\'e map} is
$T:= T_1$. 

We postulate  the following conditions for Equation (\ref{eq:per}):
\begin{description}

\item[(A1)] {\em total competition: $\pde{G_i}{x_j}\le 0, \ (i,j=1,\dots,n)$}

\item[(A2)] {\em strong self-competition: $ \sum_{k\in \mbf I (x)}
  \pde{G_k}{x_k} (t,x)<0$}

\item[(A3)]  {\em decrease of large  population:  $G_i (t,
  x) <0$ for $x_i$   sufficiently large.}

This implies existence of a {\em global attractor} for the Poincar'e map $T$.

\item[(A4)] {\em increase of small populations: $G_i(t, 0) >0$.}

\end{description}  

Under these assumptions there are two obvious candidates for a
carrying simplex for $T$, namely $\p_\K B$ and $\p_\K\Gam$,   the respective
boundaries in $\K$ of $B(0)$ and $\Gam$. 
Existence of a unique carrying simplex implies $\p_\K B=\p_\K\Gam$.

\begin{theorem}		\mylabel{th:per}
Assume system  (\ref{eq:per}) has properties {\em (A1)---(A4)}.  Then
there is a unique carrying simplex, and it is the global attractor for
the dynamics in $\K\setminus \{0\}$.
\end{theorem}
The proof, which will be given elsewhere, uses a subtle dynamical
consequence of competition discovered by Wang \& Jiang \cite[page 630]
{WangJifa}:  If $u(t), v(t)$ are solutions to Equation (\ref{eq:per})
such that for all $i$
\[
  u_i(t) < v_i(t), \qquad (s< t< s_1),
\]
then 
\[
  \ode {~} t \left(\frac{u_i}{v_i}\right)  >0, \qquad (s<t< s_1) 
\]

\begin{example}		\mylabel{th:vlper}
A competitive, periodic Volterra-Lotka  system in $\K$ of the form
\[
  \dot u_i=u_i\big(B_i(t)  -\sum_j A_{ij}(t) u_i\big),\qquad
  B_i, A_{ij} >0 
\]
satisfies (A1)---(A4)  and thus the conclusion of Theorem
\ref{th:per}. 
\end{example}

\begin{example}		\mylabel{th:vl}
Several mathematicians have investigated carrying simplex dynamics  for
competitive, autonomous Volterra-Lotka systems in $\K$  having the form
\begin{equation}		\label{eq:vl}
  \dot u_i=u_i \big(B_i -\sum_j A_{ij} u_i\big):=u_i H_i
(u_1,\dots,u_n),\quad B_i, A_{ij} >0
\end{equation} 
The best results are for $n=3$: the interesting dynamics is on a
2-dimensional cell, therefore the Poincar\'e-Bendixson theorem
\cite{Hartman} precludes any kind of chaos and makes the dynamics easy
to analyze.  The dynamics for generic systems were classified by
M.L. Zeeman \cite{Zeeman93}, with computer graphics exhibited in
Zeeman \cite{ZeemanPics}.  She proved that in many cases simple
algebraic criteria on the coefficients determine the existence of
limit cycles and Hopf bifurcations.

Van den Driessche and Zeeman \cite{vandenDriesscheZeeman04} applied
Zeeman's classification to model two competing species with species 1,
but not species 2, susceptible to disease.  They showed that if
species 1 can drive species 2 to extinction in the absence of disease,
then the introduction of disease can weaken species 1 sufficiently to
permit stable or oscillatory coexistence of both species.

 Zeeman \& Zeeman \cite {ZeemanZeeman02} showed that generically, but
not in all cases, the carrying simplex is uniquely determined by the
dynamics in the $2$-dimensional facets of $\K$.  Systems with two and
three limit cycles have been found by Hofbauer \& So \cite{HofbauerSo},
Lu \& Luo \cite{LuLuo02}), and Gyllenberg {\em et al.\
}\cite{GyllenbergYanWang06}.  No examples of Equation (\ref{eq:vl})
with four limit cycles are known.
  
More information on the dynamics of  Equation (\ref{eq:vl}) can be
found in \cite{vandenDriesscheZeeman98, 
XiaoLi00,
ZeemanE02, 
ZeemanZeeman94, 
ZeemanZeeman03}.  
\end{example}

\section*{Background}
In an important paper on competitive maps, Smith \cite{Smith86}
investigated $C^2$ diffeomorphisms $T$ of $\K$.  Under assumptions
similar to (C0)---(C5) he proved $T$ is strictly retrotone and
established the existence of the global attractor $\Gam$ and the
repulsion basin $B(0)$.  He showed that $\p_\K B(0)$ and $\p_\K\Gam$
are compact unordered invariant sets homeomorphic to the unit simplex,
and each of them contains all periodic orbits except the origin.  His
conjecture that $\p_\K B=\p_\K\Gam$ remains unproved from his
hypotheses.  He also showed that for certain types of competitive
planar maps every bounded trajectory converges, extending earlier
results of Hale \& Somolinos \cite{HaleSomolinos}, de Mottoni \&
Schiaffino \cite{deMottoniSchiaffino81}.

Using Smith's results and those of 
Hess \& Pol\'a\v cik \cite {HessPolacik}, Wang \& Jiang \cite
{WangJifa} obtained unique carrying simplices for competitive $C^2$
maps.

  
For further results on the smoothness, geometry and dynamics of
carrying simplices, see \cite {Benaim97, Mierczynski94a,
Mierczynski99, Mierczynski99a, Mierczynski99b}.


\paragraph{Mea culpa} 
Uniqueness of the carrying simplex for Equation (\ref{eq:vl}) was
claimed in Hirsch \cite{Hirsch88a}, but M.L. Zeeman \cite{ZeemanGap}
discovered an error in the proof of Proposition 2.3(d).


\end{document}